\documentclass[12pt,reqno]{amsart}
\usepackage{amssymb}
\usepackage[latin1]{inputenc}
\usepackage{amssymb}
\usepackage[all]{xy}
\usepackage{amsfonts}
\usepackage{float}
\usepackage{dsfont}
\usepackage[frenchb,english]{babel}
\usepackage{enumerate}
\usepackage{amsfonts, amsmath, amssymb, stmaryrd, latexsym}
\usepackage{icomma} 
\usepackage{pstricks}
\usepackage[all]{xy}
\usepackage{float}
\usepackage{amscd}
\usepackage{longtable}
\usepackage{geometry}
\usepackage{array}
\usepackage[frenchb,english]{babel}
\begin{document}
\theoremstyle{plain}
\newtheorem{thm}{Theorem}
\newtheorem{lem}[thm]{Lemma}
\newtheorem{prop}[thm]{Proposition}
\newtheorem{cor}[thm]{Corollary}
\theoremstyle{remark}
\newtheorem{rem}{Remark}
\newtheorem{Def}{Definition}
\newtheorem{ex}{Example}
\def\K{\mathbb{Q}(\sqrt{d},i)}
\def\NN{\mathds{N}}
\def\RR{\mathbb{R}}
\def\HH{I\!\! H}
\def\QQ{\mathbb{Q}}
\def\CC{\mathbb{C}}
\def\ZZ{\mathbb{Z}}
\def\OO{\mathcal{O}}
\def\kk{\mathbf{k}}
\def\KK{\mathbf{K}}
\def\ho{\mathcal{H}_0^{\frac{h(d)}{2}}}
\def\LL{\mathbb{L}}
\def\h{\mathcal{H}_1}
\def\hh{\mathcal{H}_2}
\def\hhh{\mathcal{H}_3}
\def\hhhh{\mathcal{H}_4}
\def\jj{\overline{\mathcal{H}_3}}
\def\jjj{\overline{\mathcal{H}_4}}
\def\L{\mathbf{k}_2^{(2)}}
\def\M{\mathbf{k}_2^{(1)}}
\def\k{\mathbf{k}^{(*)}}
\def\l{\mathds{L}}

\title[On some metabelian  $2$-group...]{On some metabelian  $2$-group whose\\ abelianization is of type $(2, 2, 2)$ and applications}
\author{A. Azizi}
\author{A. Zekhnini}
\address{Abdelmalek Azizi and Abdelkader Zekhnini, D\'epartement de Math\'ematiques\\
Facult\'e des Sciences\\
Universit\'e Mohammed 1\\
Oujda\\
Morocco}
\email{abdelmalekazizi@yahoo.fr}
\email{zekha1@yahoo.fr}
\author{M. Taous}
\address{Mohammed Taous, D\'epartement de Math\'ematiques, Facult\'e des Sciences et Techniques, Universit\'e Moulay Ismail, Errachidia, Morocco.}
\email{taousm@hotmail.com}



\maketitle
\selectlanguage{english}
\begin{abstract}
Let $G$ be some metabelian $2$-group satisfying the condition $G/G'\simeq \ZZ/2\ZZ\times\ZZ/2\ZZ\times\ZZ/2\ZZ$. In this paper,
we construct all the subgroups of $G$ of index $2$ or $4$,  we give the abelianization types  of these subgroups and we compute the kernel of the transfer map. Then we apply these  results to study the capitulation problem of the  $2$-ideal classes of some fields $\mathbf{k}$ satisfying the condition  $\mathrm{G}al(\mathbf{k}_2^{(2)}/\mathbf{k})\simeq G$, where $\mathbf{k}_2^{(2)}$ is  the second Hilbert $2$-class field of $\mathbf{k}$.
\end{abstract}
\section{Introduction}
Let $G$ be a group and let   $G'=[G, G]$ be its     \textit{derived group}, that is the subgroup of $G$ generated by the commutators $[x, y]=x^{-1}y^{-1}xy$, where $x$ and $y$ are in $G$. Let $\gamma_i(G)$ be
the i-th term of the lower central series of $G$ defined inductively  by
$\gamma_1(G)=G$ and $\gamma_{i+1}(G)=[\gamma_i(G), G]$. The group  $G$ is said to be
 \textit{nilpotent} if there exists a positive integer  $c$ such that
$\gamma_{c+1}(G)=1$; the smallest integer $c$ satisfying this equality  is called \textit{the nilpotency class} of $G$; if   the order of $G$ is $p^h$, where $p$ is a prime, then the number  $cc(G)=h-c$ is called the coclass of $G$.  It is well known that a $p$-group $G$ is nilpotent.  Recall that a group  $G$ is said to be
\textit{metabelian} if its derived group  $G'$ is abelian.  Finally, a subgroup $H$ of a  group $G$,   not reduced to an element,  is called maximal if it is the unique subgroup of $G$ distinct from  $G$  containing $H$.

Consider the  group family defined, for any  integer $n\geq1$,  as follows
\begin{align}
\begin{aligned}\label{1}
G_n=\langle \sigma,\ \tau,\ \rho: \quad & \rho^4=\sigma^8=\tau^{2^{n+2}}=1,\ \sigma^4=\tau^{2^{n+1}},\  \rho^2=\tau^{2^n}\sigma^2,\\
                                    & [\tau, \sigma]=1,\  [\rho, \sigma]=\sigma^{-2},\  [\rho, \tau]=\tau^2 \rangle.
\end{aligned}
\end{align}
 In this paper, we construct all the subgroups of $G_n$ of index $2$ or $4$,  we give the abelianization types  of these subgroups and we compute the kernel of the transfer map ${\rm V}_{G\rightarrow H}: G_n/G'_n\rightarrow H/H'$, for any subgroup $H$ of $G_n$, defined by the Artin map. Then we apply these  results to study the capitulation problem of the  $2$-ideal classes of some fields $\mathbf{k}$ satisfying the condition  $\mathrm{G}al(\mathbf{k}_2^{(2)}/\mathbf{k})\simeq G_n$, where $\mathbf{k}_2^{(2)}$ is  the second Hilbert $2$-class field of $\mathbf{k}$. Finally,  we  illustrate  our results by an example about an
imaginary bicyclic biquadratic field.
\section{ Results\label{JNT1}}
Recall first that if $x$,  $y$ and $z$ are elements of some group $G$, then  we easily  show that:
\begin{center}
$[xy, z] = [x, z]^y[y, z]$ and\\
$[x ,yz] = [x, z][x, y]^z$,\\ where $x^y=y^{-1}xy$.
\end{center}

Let $G_n$ be the  group family defined by the Equality \eqref{1}. Since $[\tau, \sigma]=1$,   $[\rho, \sigma]=\sigma^{-2}$ and   $[\rho, \tau]=\tau^2$, so $G_n'=\langle\sigma^2, \tau^2\rangle$, which is abelian.  Then  $G_n$ is metabelian and  $G_n/G_n'\simeq (2, 2, 2)$, since $\rho^2=\tau^{2^n}\sigma^2$.  Hence $G_n$ admits seven  subgroups of index $2$, denote them by $H_{i,2}$, and seven  subgroups of index $4$, denote them by $H_{i,4}$, where $1\leq i\leq7$. These subgroups,
 their derived groups and their abelianizations   are given in  the Tables \ref{26} and \ref{27}.
 \begin{table}[H]
\caption{\label{26} Subgroups of $G_n$ of index $2$.}
\begin{tabular}{c c  c  c  c }
\hline
 $i$ & \text{ conditions }& $H_{i, 2}$  & $H_{i, 2}'$ &  $H_{i, 2}/H_{i, 2}'$\\
\hline
$1$ & & $\langle\sigma, \tau \rangle$ &  $\langle1\rangle$  &  $\left(2^2, 2^{n+2}\right)$ \\ \hline
  &$n=1$ & $\langle\sigma, \rho \rangle$   & $\langle\sigma^2 \rangle$  & $\left(2, 4\right)$  \\[-1ex]
 \raisebox{2ex}{$2$}
 & $n\geq2$ & $\langle\sigma, \tau^2, \rho\rangle$   & $\langle\sigma^2, \tau^4\rangle$  & $\left(2, 2, 2\right)$  \\[1ex] \hline
 $3$ & & $\langle\rho, \tau \rangle$ & $\langle\tau^2\rangle$ &  $\left(2, 4\right)$\\ \hline
  &$n=1$ & $\langle\sigma\tau, \rho, \sigma^2 \rangle$   & $\langle(\sigma\tau)^2, \sigma^4 \rangle$  & $\left(2, 2, 2\right)$  \\[-1ex]
 \raisebox{2ex}{$4$}
 & $n\geq2$ & $\langle\sigma\tau, \rho\rangle$   & $\langle(\sigma\tau)^2\rangle$  & $\left(2, 4\right)$  \\[1ex] \hline

$5$ && $\langle\sigma\rho, \tau \rangle$ & $\langle\tau^2\rangle$ & $\left(2, 4\right)$ \\\hline

&$n=1$ & $\langle\tau\rho, \sigma \rangle$ & $\langle\sigma^2 \rangle$ & $\left(2, 4\right)$  \\[-1ex]
 \raisebox{2ex}{$6$}
& $n\geq2$ & $\langle\tau\rho, \sigma, \tau^2\rangle$ & $\langle\sigma^2, \tau^4\rangle$ & $\left(2, 2, 2\right)$  \\[1ex]\hline
&$n=1$ & $\langle\sigma\tau, \tau\rho, \sigma^2\rangle$ & $\langle\sigma^4, (\sigma\tau)^4 \rangle$ & $\left(2, 2, 2\right)$  \\[-1ex]
 \raisebox{2ex}{$7$}
& $n\geq2$ & $\langle\sigma\tau, \tau\rho\rangle$ & $\langle(\sigma\tau)^2\rangle$ & $\left(2, 4\right)$  \\[1ex]\hline
\end{tabular}
\end{table}

 \begin{table}[H]
\caption{\label{27} Subgroups of $G_n$ of index $4$.}
\begin{tabular}{c c  c  c}
 \hline
$i$ &  $H_{i, 4}$  & $H_{i, 4}'$ &  $H_{i, 4}/H_{i, 4}'$ \\ \hline
$1$  & $\langle\sigma, \tau^2 \rangle$ &  $\langle1\rangle$  &  $\left(2^{\min(n, 2)},  2^{\max(n+1, 3)}\right)$ \\ \hline
$2$  & $\langle\tau, \sigma^2 \rangle$   & $\langle1\rangle$  & $\left(2,  2^{n+2}\right)$  \\ \hline
$3$  & $\langle\rho, \tau^2 \rangle$ & $\langle\tau^4\rangle$ &  $\left(2, 4\right)$\\ \hline
&     & & $\left(4, 4\right)$ \text{ if }$n=1$ \\[-1ex]
 \raisebox{2ex}{$4$} & \raisebox{2ex}{$\langle\sigma\tau, \sigma^2 \rangle$} & \raisebox{2ex}{$\langle1\rangle$ }
 & $\left(2,  2^{n+2}\right)$ \text{ if } $n\geq2$ \\[1ex]\hline
$5$ &  $\langle\sigma\rho,\tau^2\rangle$   & $\langle\tau^4\rangle$  & $\left(2, 4\right)$  \\ \hline
$6$ &  $\langle\tau\rho,\tau^2\rangle$   & $\langle\tau^4\rangle$  & $\left(2, 4\right)$  \\ \hline
$7$ &  $\langle\sigma\tau\rho,\tau^2\rangle$   & $\langle\tau^4\rangle$  & $\left(2, 4\right)$  \\ \hline
\end{tabular}
\end{table}

To check the Tables entries we need the following lemma.
\begin{lem}\label{2}
Let $G_n=\langle\sigma, \tau, \rho\rangle$ denote the group defined above, then
\begin{enumerate}[\rm\indent 1.]
  \item $\rho^{-1}\sigma\rho= \sigma^3.$
  \item $\rho^{-1}\tau\rho=\tau^{-1}$.
  \item $\rho^2$ commutes with $\sigma$ and $\tau$.
  \item $(\tau\rho)^2=\rho^2$.
  \item $(\sigma\tau\rho)^2=(\sigma\rho)^2=\rho^2\sigma^{4}$.
  \item For all $r\in\NN$, $[\rho,\ \tau^{2^r}]=\tau^{2^{r+1}}$ and  $[\rho,\ \sigma^{2^r}]=\sigma^{-2^{r+1}}.$
\end{enumerate}
\end{lem}
\begin{proof}
1. and 2. are obvious, since $[\rho,\ \sigma]=\sigma^{-2}$ and $[\rho,\ \tau]=\tau^2$.\\
3. As $\rho^2=\tau^{2^n}\sigma^2$, so $\rho^2\in\langle\tau,\ \sigma\rangle$, which is an abelian group, because $[\tau,\ \sigma]=1$. Hence the result.\\
4. $(\tau\rho)^2=\tau\rho\tau\rho=\tau\rho^2\rho^{-1}\tau\rho=\tau\rho^2\tau^{-1}=\rho^2$.\\
5. We proceed as in 4. to prove these results.\\
6.
Since  $[\rho,\ \tau]=\tau^{2}$,
so  $[\rho,\ \tau^2]=\tau^{4}$.
By induction, we show that for all $r\in\NN^*$,   $[\rho,\ \tau^{2^r}]=\tau^{2^{r+1}}$. Similarly, we prove  that $[\rho,\ \sigma^{2^r}]=
                              \sigma^{2^{r+1}}.$
\end{proof}

Let us now prove some entries of the Tables \ref{26} and \ref{27}, using the Lemma \ref{2}.\\
$\bullet$ For $H_{1, 2}=\langle\sigma,\ \tau, G'_n\rangle=\langle\sigma,\ \tau\rangle$, we have $H_{1, 2}'=\langle1\rangle$, since $[\sigma,\tau]=1$. As $\sigma^8=\tau^{2^{n+2}}=1$ and  $\sigma^4=\tau^{2^{n+1}}$, so $H_{1, 2}/H_{1, 2}'\simeq \left(2^{\min(2, n+1)}, 2^{\max(3, n+2)}\right)=\left(2^2, 2^{n+2}\right)$.\\
$\bullet$ For $H_{2, 2}=\langle\sigma,\ \rho, G'_n\rangle=\langle\sigma,\ \rho, \ \tau^2,\ \sigma^2\rangle=\langle\sigma,\ \rho,\ \tau^2\rangle$. If $n=1$, then $\rho^2=\tau^{2}\sigma^2$,  hence $H_{2, 2}=\langle\sigma,\ \rho\rangle$. Therefore $H_{2, 2}'=\langle\sigma^2\rangle$ and  $H_{2, 2}/H_{2, 2}'\simeq(2, 4)$.  If $n\geq2$, then  $H_{2, 2}=\langle\sigma,\ \rho,\ \tau^2\rangle$. Therefore $H_{2, 2}'=\langle\sigma^2,\ \tau^4\rangle$ and  $H_{2, 2}/H_{2, 2}'\simeq(2, 2, 2)$. \\
$\bullet$ For $H_{1, 4}=\langle\sigma,\  G'_n\rangle=\langle\sigma,\ \sigma^2,\ \tau^2\rangle=\langle\sigma,\  \tau^2\rangle$, we have $H_{1, 4}'=\langle1\rangle$, since $[\sigma,\tau]=1$. As $\sigma^8=\tau^{2^{n+2}}=1$ and  $\sigma^4=\tau^{2^{n+1}}$, so $H_{1, 4}/H_{1, 4}'\simeq \left(2^{\min(2, n)}, 2^{\max(3, n+1)}\right)$.\\
$\bullet$ For $H_{2, 4}=\langle\tau,\  G'_n\rangle=\langle\tau,\  \tau^2,\ \sigma^2\rangle=\langle\tau,\ \sigma^2\rangle$, we have   $H_{2, 4}'=\langle1\rangle$, hence  $H_{2, 4}/H_{2, 4}'\simeq\left(2^{\min(1, n+1)}, 2^{\max(2 n+2)}\right)=\left(2, 2^{n+2}\right)$. \\
The other entries are checked similarly.
\begin{prop}\label{7}
Let $G_n$ be the  group family defined by the Equality \eqref{1}, then
\begin{enumerate}[\rm\indent 1.]
\item The order of $G_n$ is $2^{n+5}$ and that of $G_n'$ is $2^{n+2}$.
\item The coclass of $G_n$ is $3$ and its nilpotency class is $n+2$.
\end{enumerate}
\end{prop}
\begin{proof}
1. Since $\sigma^8=\tau^{2^{n+2}}=1$ and  $\sigma^4=\tau^{2^{n+1}}$, then
 $$\langle\sigma,\ \tau \rangle\simeq\left(2^{\min(2, n+1)}, 2^{\max(3,  n+2)}\right)=\left(2^2, 2^{n+2}\right).$$
   Moreover, as $\rho^2=\tau^{2^n}\sigma^2$  and $\rho^4=\sigma^8=\tau^{2^{n+2}}=1$, so
   $$\langle\sigma,\ \tau,\ \rho \rangle\simeq\left(2^2,\ 2^{n+2},\ 2\right).$$
   Thus $|G_n|=2^{n+5}$. Similarly, we prove that $|G_n'|=2^{n+2}$, since  $G_n'=\langle\sigma^2,\ \tau^2\rangle$.

2.  The lower central series of $G_n$ is defined inductively by $\gamma_1(G_n)=G_n$ and $\gamma_{i+1}(G_n)=[\gamma_i(G_n),G_n]$, that is the subgroup of $G_n$ generated by the set $\{[a, b]=a^{-1}b^{-1}ab/ a\in \gamma_i(G_n), b\in G_n\}$, so  the coclass of $G_n$ is defined to be $cc(G_n) = h-c$, where $|G_n|=2^h$ and  $c=c(G_n)$ is the nilpotency class of $G_n$. We easily get \\
$\gamma_1(G)=G$.\\
$\gamma_2(G)=G'=\langle\sigma^2, \tau^2\rangle$.\\
$\gamma_3(G)=[G',G]=\langle\sigma^4, \tau^4\rangle$.\\
 Then   Lemma \ref{2}(6) yields that $\gamma_{j+1}(G)=[\gamma_j(G),G]=\langle\sigma^{2^{j}}, \tau^{2^{j}}\rangle$. Hence, if  we put $\upsilon=\max(n+1, 2)=n+1$, then we get $\gamma_{\upsilon+2}(G)=\langle\sigma^{2^{\upsilon+1}}, \tau^{2^{\upsilon+1}}\rangle=\langle1\rangle$ and $\gamma_{\upsilon+1}(G)=\langle\sigma^{2^{\upsilon}}, \tau^{2^{\upsilon}}\rangle\neq\langle1\rangle$. As $\mid G\mid=2^{n+5}$, so $$c(G)=n+2\text{ and } cc(G) = n+5-(n+2)=3.$$
 \end{proof}
 We contenue with the following results.
 \begin{prop}[\cite{Mi89}]
 Let $H$ be a normal subgroup of $G_n$. For $g\in G_n$, put  $f=[\langle
g\rangle.H:H]$ and let   $\{x_1,  x_2,
\ldots, x_t\}$ be a set of representatives of $G/\langle
g\rangle H$.  The transfer map ${\rm V}_{G\rightarrow H}: G/G'\rightarrow H/H'$ is given by the following formula
\begin{equation}\label{JNT7}
{\rm V}_{G\rightarrow H}(gG')=\prod_{i=1}^tx_i^{-1}g^fx_i.H'.
\end{equation}
 \end{prop}
 Easily, we prove the following corollaries.
 \begin{cor}\label{3}
  Let $H$ be a  subgroup of $G_n$ of index $2$. If  $G_n/H=\{1, zH\}$, then \
${\rm V}_{G\rightarrow H}(gG'_n)=\left\{
\begin{array}{ll}
gz^{-1}gz.H'=g^2[g,\ z].H' & \text{ if }g\in H,\\
g^2.H' & \text{ if }g\not\in H.
\end{array}
\right.$
 \end{cor}
  \begin{cor}\label{4}
  Let $H$ be a normal  subgroup of $G_n$ of index $4$. If  $G_n/H=\{1, zH, z^2H, z^3H\}$, then
 $${\rm V}_{G\rightarrow H}(gG'_n)=\left\{
\begin{array}{ll}
gz^{-1}gz^{-1}gz^{-1}gz^3.H' & \text{ if }g\in H,\\
g^4.H' & \text{ if }gH=zH,\\
g^2z^{-1}g^2z.H' & \text{ if }g\not\in H \text{ and }gH\neq zH.
\end{array}
\right.$$
 \end{cor}
   \begin{cor}\label{5}
  Let $H$ be a normal  subgroup of $G_n$ of index $4$. If  $G_n/H=\{1, z_1H, z_2H, z_3H\}$ with $z_3=z_1z_2$, then \
$${\rm V}_{G\rightarrow H}(gG'_n)=\left\{
\begin{array}{ll}
gz_1^{-1}gz_1z_2^{-1}gz_1^{-1}gz_1z_2.H' & \text{ if }g\in H,\\
g^2z_i^{-1}g^2z_i.H' & \text{ if } gH=z_jH \text{ with }i\neq j.
\end{array}
\right.$$
 \end{cor}
\begin{thm}\label{6}
Keep the previous notations.  Then  we have:
\begin{enumerate}[\rm\indent 1.]
\item $\ker\mathrm{V}_{G_n \rightarrow H_{1, 2}}=\langle\tau G'_n\rangle$.
\item $\ker \mathrm{V}_{G_n \rightarrow H_{2, 2}}=
\left\{
\begin{array}{ll}
\langle\sigma G_n', \rho G'_n\rangle & \text{ if }n=1.\\
\langle\sigma G_n', \tau\rho G'_n\rangle & \text{ if }n\geq2.
\end{array}
\right.$
\item $\ker \mathrm{V}_{G_n \rightarrow H_{3, 2}}=\langle\tau G'_n,  \rho G'_n\rangle$.
\item $\ker \mathrm{V}_{G_n \rightarrow H_{4, 2}}=
\left\{
\begin{array}{ll}
\langle\tau\rho G'_n, \sigma\rho G'_n\rangle & \text{ if }n=1.\\
\langle\sigma\tau G'_n, \rho G'_n\rangle & \text{ if }n\geq2.
\end{array}
\right.$
\item $\ker \mathrm{V}_{G_n \rightarrow H_{5, 2}}=\langle\tau G'_n, \sigma\rho G'_n\rangle$.
\item $\ker \mathrm{V}_{G_n \rightarrow H_{6, 2}}=
\left\{
\begin{array}{ll}
\langle\sigma G'_n, \tau\rho G'_n\rangle & \text{ if }n=1.\\
\langle\rho G'_n, \sigma G'_n\rangle & \text{ if }n\geq2.
\end{array}
\right.$
\item $\ker \mathrm{V}_{G_n \rightarrow H_{7, 2}}=
\left\{
\begin{array}{ll}
\langle\sigma\tau G'_n, \rho G'_n\rangle & \text{ if }n=1.\\
\langle\sigma\tau G_n', \tau\rho G'_n\rangle & \text{ if }n\geq2.
\end{array}
\right.$
\item  For all $1\leq i\leq7$, $ \ker \mathrm{V}_{G_n \rightarrow H_{i, 4}}=G_n/G'_n$.
\end{enumerate}
\end{thm}
\begin{proof}  We prove only some assertions  of the theorem.\\
1. We know, from the Table \ref{26}, that $H_{1, 2}=\langle\sigma, \tau \rangle$, then $G_n/H_{1, 2}=\{1,\ \rho H_{1, 2}\}$ and $H_{1, 2}'=\langle1\rangle$. Hence, by Corollary \ref{3} and Lemma \ref{2}, we get
\begin{enumerate}[\indent $\ast$]
\item $\mathrm{V}_{G_n \rightarrow H_{1, 2}}(\sigma G'_n)=\sigma^2[\sigma,\ \rho]H_{1, 2}'=\sigma^2\sigma^{2}H_{1, 2}'=\sigma^{4}H_{1, 2}'$.
\item $\mathrm{V}_{G_n \rightarrow H_{1, 2}}(\tau G'_n)=\tau^2[\tau,\ \rho]H_{1, 2}'= \tau^2\tau^{-2}H_{1, 2}'=H_{1, 2}'$.
 \item $\mathrm{V}_{G_n \rightarrow H_{1, 2}}(\rho G'_n)=\rho^2H_{1, 2}'\neq H_{1, 2}'$.
 \item $\mathrm{V}_{G_n \rightarrow H_{1, 2}}(\sigma\tau G'_n)=(\sigma\tau)^2[\sigma\tau,\ \rho]H_{1, 2}'=\tau^{4}H_{1, 2}'$.
\item $\mathrm{V}_{G_n \rightarrow H_{1, 2}}(\tau\rho G'_n)=(\tau\rho)^2H_{1, 2}'= \rho^2H_{1, 2}'\neq H_{1, 2}'$.
 \item $\mathrm{V}_{G_n \rightarrow H_{1, 2}}(\sigma\rho G'_n)=(\sigma\rho)^2H_{1, 2}'=\rho^2\sigma^4 H_{1, 2}'\neq H_{1, 2}'$.
 \item $\mathrm{V}_{G_n \rightarrow H_{1, 2}}(\sigma\tau\rho G'_n)=(\sigma\tau\rho)^2H_{1, 2}'=\rho^2\sigma^4 H_{1, 2}'\neq H_{1, 2}'$.
 \end{enumerate}
 Therefore $\ker\mathrm{V}_{G_n \rightarrow H_{1, 2}}=\langle\tau G'_n\rangle.$\\
 8. We know, from the Table \ref{27}, that $H_{1, 4}=\langle\sigma, \tau^2\rangle$, then $G_n/H_{1, 4}=\{1,\  \tau H_{1, 4},\ \rho H_{1, 4},\ \tau\rho H_{1, 4}\}$ and $H_{1, 4}'=\langle1\rangle$. Hence Corollary \ref {5} and Lemma \ref{2} yield that
\begin{enumerate}[\indent $\ast$]
\item $\mathrm{V}_{G_n \rightarrow H_{1, 4}}(\sigma G'_n)=\sigma^8H_{1, 4}'=H_{1, 4}'$.
\item $\mathrm{V}_{G_n \rightarrow H_{1, 4}}(\tau G'_n)=\tau^2\rho^{-1}\tau^2\rho H_{1, 4}'=H_{1, 4}'$.
 \item $\mathrm{V}_{G_n \rightarrow H_{1, 4}}(\rho G'_n)=\rho^4H_{1, 4}'= H_{1, 4}'$.
 \end{enumerate}
 Therefore $\ker\mathrm{V}_{G_n \rightarrow H_{1, 4}}=\langle\sigma G'_n,\ \tau G'_n,\ \rho G'_n\rangle.$\\
 The other assertions are proved similarly.
\end{proof}
\section{Applications $\label{009}$}
Let $\mathbf{k}$ be a number field and
     $\mathrm{C}_{\mathbf{k}, 2}$ be its $2$-class group, that
is the $2$-Sylow subgroup of the ideal class group $\mathrm{C}_{\mathbf{k}}$ of $\mathbf{k}$, in the wide sens. Let
$\mathbf{k}^{(1)}_2$ be the  Hilbert $2$-class field of  $\mathbf{k}$ in the wide sens.  Then the  Hilbert $2$-class field tower of $\mathbf{k}$ is defined inductively  by:
$\mathbf{k}^{(0)}_2=\mathbf{k}$ and
$\mathbf{k}^{(n+1)}_2=(\mathbf{k}^{(n)}_2)^{(1)}$, where $n$ is a positive integer. Let
$\mathbb{M}$ be an unramified  extension of $\mathbf{k}$ and
$\mathrm{C}_{\mathbb{M}}$ be the subgroup of $\mathrm{C}_{\mathbf{k}}$ associated to $\mathbb{M}$ by the class field theory. Denote by
$j_{\mathbf{k} \rightarrow \mathbb{M}}: \mathrm{C}_{\mathbf{k}} \longrightarrow \mathrm{C}_{\mathbb{M}}$  the homomorphism  that associate  to the class of an ideal  $\mathcal{A}$ of  $\mathbf{k}$ the class of the ideal generated by $\mathcal{A}$ in $\mathbb{M}$,  and by
$\mathcal{N}_{\mathbb{M}/\mathbf{k}}$ the norm of the extension
$\mathbb{M}/\mathbf{k}$.

Throughout all this section, assume that $\mathrm{G}al(\L/\kk)\simeq G_n$. Hence, according to the class field theory,  $\mathrm{C}_{\mathbf{k}, 2}\simeq G_n/G'_n\simeq (2, 2, 2)$, thus   $\mathrm{C}_{\mathbf{k},
2}=\langle\mathfrak{a},  \mathfrak{b}, \mathfrak{c}\rangle\simeq \langle \sigma G'_n,
\tau G'_n, \rho G'_n\rangle$, where  $(\mathfrak{a}, \mathbf{k}_2^{(2)}/\mathbf{k})=\sigma G'_n$, $(\mathfrak{b}, \mathbf{k}_2^{(2)}/\mathbf{k})=\tau G'_n$
and  $(\mathfrak{c}, \mathbf{k}_2^{(2)}/\mathbf{k})=\rho G'_n$, with $(\ .\ ,
\mathbf{k}_2^{(2)}/\mathbf{k})$ denotes the Artin symbol in
$\mathbf{k}_2^{(2)}/\mathbf{k}$.

 It is well known that each subgroup  $H_{i,j}$, where $1\leq i \leq 7$ and $j=2$ or $4$,  of  $\mathrm{C}_{\mathbf{k}, 2}$ is associated, by class field theory, to a unique unramified extension
 $\mathbf{K}_{i,j}$ of  $\mathbf{k}_2^{(1)}$ such that  $H_{i,j}/H'_{i,j}\simeq\mathrm{C}_{\mathbf{K}_{i,j}, 2}$.

Our goal  is  to study  the capitulation problem of the  $2$-ideal classes of $\mathbf{k}$ in its unramified quadratic and biquadratic extensions $\mathbf{K}_{i, 2}$ and $\mathbf{K}_{i, 4}$.
By the class field theory,  the kernel of  $j_{\mathbf{k} \rightarrow \mathbb{M}}$, $\ker j_{\mathbf{k} \rightarrow \mathbb{M}}$,  is determined by the kernel of the transfer map ${\rm V}_{G\rightarrow H}: G/G'\rightarrow H/H'$, where  $G=\mathrm{Gal}(\mathbf{k}_2^{(2)}/\mathbf{k})$ and  $H=\mathrm{Gal}(\mathbb{M}_2^{(2)}/\mathbb{M})$.

\begin{thm}\label{003}
Keep the previous notations.  Then  we have:
\begin{enumerate}[\rm\indent 1.]
\item $\ker j_{\mathbf{k} \rightarrow \mathbf{K}_{1, 2}}=\langle\mathfrak{b}\rangle$.
\item $\ker j_{\mathbf{k} \rightarrow \mathbf{K}_{2, 2}}=
\left\{
\begin{array}{ll}
\langle\mathfrak{a}, \mathfrak{c}\rangle & \text{ if }n=1.\\
\langle\mathfrak{a}, \mathfrak{bc}\rangle & \text{ if }n\geq2.
\end{array}
\right.$
\item $\ker j_{\mathbf{k} \rightarrow \mathbf{K}_{3, 2}}=\langle\mathfrak{b}, \mathfrak{c}\rangle$.
\item $\ker j_{\mathbf{k} \rightarrow \mathbf{K}_{4, 2}}=
\left\{
\begin{array}{ll}
\langle \mathfrak{ac}, \mathfrak{bc}\rangle & \text{ if }n=1.\\
\langle\mathfrak{ab}, \mathfrak{c}\rangle & \text{ if }n\geq2.
\end{array}
\right.$
\item $\ker j_{\mathbf{k} \rightarrow \mathbf{K}_{5, 2}}=\langle\mathfrak{b}, \mathfrak{ac}\rangle$.
\item $\ker j_{\mathbf{k} \rightarrow \mathbf{K}_{6, 2}}=
\left\{
\begin{array}{ll}
\langle\mathfrak{a}, \mathfrak{bc}\rangle & \text{ if }n=1.\\
\langle\mathfrak{a}, \mathfrak{c}\rangle & \text{ if }n\geq2.
\end{array}
\right.$
\item $\ker j_{\mathbf{k} \rightarrow \mathbf{K}_{7, 2}}=
\left\{
\begin{array}{ll}
\langle\mathfrak{ab}, \mathfrak{c}\rangle & \text{ if }n=1.\\
\langle\mathfrak{ab}, \mathfrak{bc}\rangle & \text{ if }n\geq2.
\end{array}
\right.$
\item  For all $1\leq i\leq7$, $\ker j_{\mathbf{k} \rightarrow \mathbf{K}_{i, 4}}=\mathrm{C}_{\mathbf{k}, 2}$.
\item The $2$-class group of $\M$ is of type $\left(2,\ 2^{n+1}\right)$.
\item The  Hilbert $2$-class field tower of $\mathbf{k}$ stops at $\L$.
\end{enumerate}
\end{thm}
\begin{proof} According to the Theorem \ref{6}, we have\\
(1) Since  $\ker\mathrm{V}_{G_n \rightarrow H_{1, 2}}=\langle\tau G'_n\rangle$, so $\ker j_{\mathbf{k} \rightarrow \mathbf{K}_{1, 2}}=\langle\mathfrak{b}\rangle$.\\
(2) Similarly, as
$\ker \mathrm{V}_{G_n \rightarrow H_{2, 2}}=
\left\{
\begin{array}{ll}
\langle\sigma G_n', \rho G'_n\rangle & \text{ if }n=1,\\
\langle\sigma G_n', \tau\rho G'_n\rangle & \text{ if }n\geq2,
\end{array}
\right.$\\
then $\ker j_{\mathbf{k} \rightarrow \mathbf{K}_{2, 2}}=
\left\{
\begin{array}{ll}
\langle\mathfrak{a}, \mathfrak{c}\rangle & \text{ if }n=1.\\
\langle\mathfrak{a}, \mathfrak{bc}\rangle & \text{ if }n\geq2.
\end{array}
\right.$\\
The other assertions are proved similarly.\\
(9) It is well known that $\mathrm{C}_{\M, 2}\simeq G_n'$, where $\mathrm{C}_{\M, 2}$ is the $2$-class group of $\M$. As $G_n'\simeq \langle\sigma^2,\ \tau^2\rangle\simeq \left(2^{\min(1, n)},\ 2^{\max(2, n+1)}\right)\simeq \left(2,\ 2^{n+1}\right)$, since $\sigma^4=\tau^{2^{n+1}}$; so the result derived.\\
(10) $H_{1, 4}$, $H_{2, 4}$ and $H_{4, 4}$ are the three subgroups of index $2$ of the group $H_{1, 2}$, then $\mathbf{K}_{1, 4}$, $\mathbf{K}_{2, 4}$ and $\mathbf{K}_{4, 4}$ are the three unramified quadratic extensions of $\mathbf{K}_{1, 2}$. On the other hand, the $2$-class groups of these fields are of rank $2$, since,  by the class field theory, $\mathrm{C}_{\mathbf{K}_{i, j}, 2}\simeq H_{i, j}/H_{i, j}'$ with $i=1$, $2$ or $4$ and $j=2$ or $4$. Thus the Tables \ref{26} and \ref{27} imply that $\mathrm{C}_{\mathbf{K}_{1, 2}, 2}\simeq \left(2^2,\ 2^{n+2}\right)$ and $\mathrm{C}_{\mathbf{K}_{2, 4}, 2}\simeq \left(2,\ 2^{n+2}\right)$. Hence $h_2(\mathbf{K}_{2, 4})=\frac{h_2(\mathbf{K}_{1, 2})}{2}$, where  $h_2(K)$ denotes the $2$-class number of the field $K$. Therefore  we can apply Proposition 7 of \cite{B.L.S-98}, which says that $\mathbf{K}_{1, 2}$ has an abelian  $2$-class field tower if and only if it has a quadratic unramified extension $\mathbf{K}_{2, 4}/\mathbf{K}_{1, 2}$ such that $h_2(\mathbf{K}_{2, 4})=\frac{h_2(\mathbf{K}_{1, 2})}{2}$.  Thus  $\mathbf{K}_{1, 2}$ has abelian  $2$-class field tower which terminates at the first stage; this implies that the $2$-class field tower of $\kk$ terminates at $\kk_2^{(2)}$, since  $\kk\subset \mathbf{K}_{1, 2}$.  Moreover, we know, from Proposition \ref{7}, that $\left|G_n\right|=2^{n+5}$ and $\left|G_n'\right|=2^{n+2}$,  hence $\kk_2^{(1)}\neq\kk_2^{(2)}$.
\end{proof}
\section{Example}
Let $p_1\equiv p_2\equiv 5 \pmod8$ be different  primes.  Denote by $\kk$ the imaginary bicyclic biquadratic field $\K$, where $d=2p_1p_2$. Let $\kk_2^{(1)}$ be the Hilbert 2-class field of $\kk$,  $\L$ its second Hilbert 2-class field and $G$ be the Galois group of $\kk_2^{(2)}/\kk$. According to  \cite{AzTa08},  $\kk$ has an elementary abelian 2-class group $\mathbf{C}_{\kk, 2}$ of rank 3, that is of type $(2, 2, 2$). Put $\KK=\kk(\sqrt 2)=\QQ(\sqrt2, \sqrt{p_1p_2}, \sqrt{-1})$, and let $q$ denote the unit index of $\KK^+=\QQ(\sqrt2, \sqrt{p_1p_2})$. Denote by $h_2(-p_1p_2)$ (resp. $h_2(p_1p_2)$) the $2$-class number of $\QQ(\sqrt{-p_1p_2})$ (resp. $\QQ(\sqrt{p_1p_2})$), then, from \cite{Ka76}, $h_2(-p_1p_2)=2^{m+1}$ with $m\geq2$, and $h_2(p_1p_2)=2^{n}$ with $n\geq1$.  Assume that $q=2$, then, by \cite[Lemma 6]{AZT14-2},  $m=2$ and  $n\geq1$, and by   \cite[Theorem 2]{AZT14-2},  $G\simeq G_n$. The following result is proved in \cite{AZT14-2}, and we give it here to illustrate the results  shown above. For more details, the reader can see \cite{AZT14-2}.
\begin{thm}
Let $p_1\equiv p_2\equiv 5 \pmod8$ be two different  primes.  Put  $\kk=\QQ(\sqrt{2p_1p_2},\ i)$,  $\kk$ has fourteen unramified extensions within his first  Hilbert $2$-class field, $\kk_2^{(1)}$ $($see \cite{AZT12-1}$)$. Denote by  $\mathbf{C}_{\kk, 2}$ the $2$-class group of $\kk$. Then the following assertions hold.
\begin{enumerate}[\rm\indent 1.]
  \item Exactly  four elements of $\mathbf{C}_{\kk, 2}$   capitulate in each unramified quadratic extension of $\kk$, except one where there are only two.
  \item All the $2$-classes of $\kk$  capitulate in each unramified biquadratic extension of $\kk$.
  \item The Hilbert $2$-class field tower of $\kk$ stops at $\kk_2^{(2)}$ $($see \cite{AZT-3}$)$.
  \item $\mathbf{C}_{\M, 2}\simeq (2,\ 2^{n+1})$.
  \item The  coclass of $G$ is $3$ and its nilpotency class is $n+2$.
  \item The $2$-class groups of the unramified quadratic extensions of $\kk$ are of type $(2,\ 4)$,  $(2,\ 2,\ 2)$ or $(4,\ 2^{n+2})$.
  \item The types of the $2$-class groups of the unramified biquadratic extensions of $\kk$ are exactly those given by the Table $\ref{27}$.
  \end{enumerate}
\end{thm}

\end{document}